\documentclass{amsart}

\usepackage{amsmath,amssymb,url}
\numberwithin{equation}{section}

\theoremstyle{plain}
\newtheorem{theorem}{Theorem}[section]
\newtheorem{lemma}[theorem]{Lemma}
\newtheorem{proposition}[theorem]{Proposition}
\newtheorem{corollary}[theorem]{Corollary}

\theoremstyle{definition}
\newtheorem{definition}[theorem]{Definition}
\newtheorem{example}[theorem]{Example}
\newtheorem{remark}[theorem]{Remark}

 \newcommand{\N}{{\mathbb{N}}}

\begin{document}
 
\title[Stable varieties of semigroups and groupoids]{Stable varieties of semigroups and groupoids}
\author[Sl. Shtrakov]{Slavcho Shtrakov}
 \email{shtrakov@swu.bg} \urladdr{http://shtrakov.swu.bg/}
\address{Department of Computer Science, South-West University,  2700 Blagoevgrad, Bulgaria}

\author[J. Koppitz]{J\"org  Koppitz}
\email{koppitz@uni-potsdam.de}
\address{Institute of Mathematics, University of
Potsdam, 14415 Potsdam}

\subjclass[2010]{Primary: 08B05; Secondary: 08A02, 03C05}

\keywords{ composition of terms,
  essential position in terms, stable variety}
  
\begin{abstract}
 The paper deals with $\Sigma$-composition and $\Sigma$-essential composition of terms which lead to stable and s-stable varieties of algebras.
A full description of all stable   varieties of semigroups, commutative and idempotent groupoids is obtained.
 We
use an abstract reduction system which simplifies the presentations of terms of type $\tau=(2)$  to study the variety of idempotent groupoids and  s-stable varieties of groupoids. S-stable varieties are a variation of stable varieties, used to highlight replacement of subterms  of a term in a deductive system instead of the usual replacement of variables by terms.

\end{abstract}
\maketitle

\section{Introduction}
Let ${\mathcal{F}}$ be any finite set of $ operation\ symbols. $ Let
$\tau:{\mathcal{F}}\to N$ be a
mapping into the non-negative integers; for $f\in{\mathcal{F}},$ the
number $\tau(f)$ will denote the \emph{arity } of the operation
symbol $f.$ The pair $(\mathcal{F},\tau)$ is called a \emph{type} or
\emph{signature}. If it is obvious what the set ${\mathcal{F}}$ is,
we will write ``$ type\ \tau$". The set of symbols of arity $p$ is
denoted by ${\mathcal{F}}_p.$

 Let $X=\{x_1,x_2,\ldots\}$  be a countable  set of variables, and let $\tau$
be a type with the set of operation symbols ${\mathcal
F}.$ The set $W_{\tau}(X)$ of \emph{
terms of type $\tau$ } with\ variables from $X $ is the
smallest set such that
 $X\cup\mathcal{F}_0\subseteq W_\tau(X)$ and  if $f $ is an $n$-ary operation symbol, and
$t_1,\ldots,t_{n}\in W_\tau(X)$ are terms, then  $f(t_1,\ldots,
t_{n})\in W_\tau(X)$.

 If $f\in\mathcal{F}$, then $f^\mathcal{A}$ denotes a $\tau(f)$-ary operation on the
set $A.$
 An \emph{algebra}
 ${\mathcal A} = \langle A; \mathcal{F}^{\mathcal A}\rangle$
  of type $\tau$ is a pair consisting of a set $A$ and an
  indexed  set $\mathcal{F}^{\mathcal{A}}$ of operations, defined on $A$.
   If $s, t \in W_{\tau}(X)$,
    then the pair $s \approx t$ is called an \emph{identity} of type
$\tau$ which is satisfied in the algebra  ${\mathcal A}$,  ${\mathcal A}\models t\approx s$ iff $t^A=s^A$. 

The operators $Id$ and $Mod$ are defined
     for classes of algebras $K$ and
     for sets of identities $\Sigma$ as
     follows
\begin{align*}
Id (K)&= \{t \approx s \mid {\mathcal A} \in K~  \Rightarrow  {\mathcal A} \models t \approx s\}, \text{ and}\\
Mod(\Sigma)&= \{ {\mathcal A} \mid   t \approx s \in  \Sigma~ \Rightarrow  {\mathcal A} \models t \approx s\}.
\end{align*}

The fixed points with respect to the closure operators $IdMod$ and $ModId$ are called  \emph{equational theories} and  \emph{varieties of algebras,} respectively.

In Section \ref{sec2} we introduce the    inductive, positional  and    $\Sigma$-composition of terms.
 
 We   apply the concept of
 $\Sigma$-composition of terms  to
study the stable varieties of semigroups (see Theorem \ref{t3.1}). We prove that a variety $V$ of semigroups is stable if and only if $Id(V)$ contains an identity of the form $(x_1x_2)x_3\approx x_ix_j$ with $1\leq i<j\leq 3$. We present a complete list of all the stable varieties of semigroups (see Theorem \ref{t3.2}).

An abstract Reduction System (ARS) for terms, which reduces  the   complexity of terms by such traditional  measures   as depth and length  is introduced in Section \ref{sec3}. 

The varieties of commutative and idempotent groupoids are stable which is shown in Section \ref{sec5}.

 We present stronger conditions for stability of varieties which successfully work in the variety of groupoids. These conditions allow us to define and  study the s-stable varieties of groupoids in Section \ref{sec6} (see Theorem \ref{t6.2}).

\section{Compositions of terms}\label{sec2}

If $t$ is a term, then the set $var(t)$ consisting of those elements
of $X$ which occur in $t$ is called the set of
 \emph{input variables (or variables)}  in $t$.  If
$t=f(t_1,\ldots,t_n)$ is a non-variable term, then $f$ is the \emph{
root symbol (root)} of $t$.

 For a term  $t\in W_\tau(X)$ the set   $Sub(t)$ of its
subterms
 is defined as follows:
if $t\in X\cup \mathcal{F}_0$, then $Sub(t)=\{t\}$ and if
$t=f(t_1,\ldots,t_{n})$, then $Sub(t)=\{t\}\cup
Sub(t_1)\cup\ldots\cup Sub(t_n).$

 Let $r,s,t\in W_\tau(X)$ be   terms of
type $\tau$. By $t(r \leftarrow s)$ we   denote the term,
obtained by simultaneous replacement of  every occurrence of $r $ as
a
 subterm of $t$ by $s$.  This term is called
 the
\emph{inductive composition }\cite{sht1} of the terms $t$ and $r$, by $s$. 
If $r_i\notin Sub(r_j)$ when $i\neq j$, then $t(r_1\leftarrow s_1,\ldots,r_m\leftarrow s_m)$ means the inductive composition of
$t,r_1,\ldots,r_m$ by $s_1,\ldots,s_m$, respectively. In the particular case when
$r_j=x_j$ for $j=1,\ldots,m$ and $var(t)=\{x_1,\ldots,x_m\}$ we will
briefly write $t(s_1,\ldots,s_m)$ instead of $t(x_1\leftarrow
s_1,\ldots,x_m\leftarrow s_m)$.

Any term can be regarded  as a tree with nodes labeled as the
operation symbols and leaves   labeled as variables or nullary
operation symbols (see Figure \ref{f-EssPos}, below).

 Let $\tau$ be a type and $\mathcal{F}$ be its set of operation
symbols. Denote by  
\[\N_\tau=\{m\in \N \mid m\leq \max_{f\in \mathcal{F}}\tau(f)\}.\]  Let
$\N_\tau^*$
be the set of all finite strings over $\N_\tau.$ The set
$\N_\tau^*$ is naturally ordered by $p\preceq q \iff p$\ is a
prefix of $q.$ The Greek letter $\varepsilon$, as usual denotes the
empty word (string) over $\N_\tau.$

To distinguish between different occurrences of the same operation
symbol in a term $t$ we assign to each occurrence of an operation symbol a
position. Usually positions are finite
sequences (strings) over $\N_\tau.$ Each position is assigned to
a node of the tree diagram of $t$, starting with the empty sequence
$\varepsilon$ for the root and using the integer $j$, $1\leq j\leq
n$ for the $j$-th branch of an $n$-ary operational symbol $f$.
Inductively, let the position $p=a_1a_2\ldots a_s\in \N^*_\tau$ be
assigned to a node of $t$ labeled by the $n$-ary operational
symbol $f$. Then the position assigned to the $j$-th child of this
node is $a_1a_2\ldots a_sj$. The set of positions
of a term $t$ is
denoted by $Pos(t)$.

  Let $t\in W_\tau(X)$ be a term of type $\tau$ and let $sub_t:Pos(t)\to
Sub(t)$  be the function  which maps each position in a term $t$ to
the subterm of $t$, whose root node occurs at that position.

 Let $t,r\in W_\tau(X)$ be two terms of type
$\tau$ and let $p\in Pos(t)$ be a position in $t.$  The \emph{positional
composition} \cite{sht1} of $t$ and $r$ on $p$ is the term $s=t(p;r)$ obtained
from $t$ by replacing   the term $sub_t(p)$ by $r$  on the position
$p$, only.
We will use notation $t(p,q;r)$ for the composition $t(p;r)(q;r)$ when
$p\not\preceq q\  \&\  q\not\preceq
p$ and more generally, if  $S=\langle p_1,\ldots,p_m\rangle\in Pos(t)^m$ with
$(\forall\ 
 i,j\leq m)\    (i\neq j \Rightarrow p_i\not\preceq p_j\  \&\  p_j\not\preceq
p_i)$
 then $t(S;r)=t(p_1,\ldots,p_m;r)= t(p_1;r)\ldots (p_m;r)$. If $T=\langle
t_1,\ldots,t_m\rangle\in W_\tau(X)^m$ then $t(S;T)= 
t(p_1;t_1)\ldots (p_m;t_m)$.

Let  $X_n=\{x_1,\ldots,x_n\}$  be a finite set of variables in $X$. Then we denote by $W_\tau(X_n)$ the set  $W_\tau(X_n)=\{t\in W_\tau(X)\ |\ var(t)\subseteq X_n\}$ of terms.

 Let $\Sigma\subseteq Id(\tau)$, $t\in
W_\tau(X_n)$ be an $n$-ary term of type $\tau$, $\mathcal{A}=\langle
A,\mathcal{F}\rangle$ be an algebra of type $\tau$ and let $x_i\in
var(t)$ be a variable which occurs in $t.$
  The
variable $x_i$ is called \emph{essential} \cite{sht}   in $t$ with respect to
 the algebra $\mathcal{A}$ if there are $n+1$ elements
$a_1,\ldots,a_{i-1},a,b,a_{i+1},\ldots,a_n\in A$  such that
\[t^\mathcal{A}(a_1,\ldots,a_{i-1},a,a_{i+1},\ldots,a_n)\neq
t^\mathcal{A}(a_1,\ldots,a_{i-1},b,a_{i+1},\ldots,a_n).\] The set of
all essential variables in $t$ with respect to  $\mathcal{A}$ is denoted by $Ess(t,\mathcal{A})$. $Fic(t,\mathcal{A})$ denotes the
set of all variables in $var(t)$, which are not essential with
respect to $\mathcal{A}$, called fictive variables.

Let $\Sigma$ be a set of identities of type $\tau$. Then ${\mathcal A}  \models
\Sigma$
means that ${\mathcal A}  \models t \approx s$ for all $ t \approx
s\in\Sigma$. For $t,s\in W_\tau(X)$ we say $\Sigma$ yields  $t\approx s$ (write:
$\Sigma\models t\approx s$)  if, given any algebra $\mathcal{A}$,\ 
$\mathcal{A}\models \Sigma\ \Rightarrow\ \mathcal{A}\models t\approx s.$

A variable $x_i$ is said to be \emph{
$\Sigma$-essential} \cite{sht1}  in a term $t$  if there is an algebra
$\mathcal{A}$, such that $\mathcal{A}\models \Sigma$ and $ x_i\in
Ess(t,\mathcal{A}).$ The set of all $\Sigma$-essential variables in
$t$ is denoted by $Ess(t,\Sigma).$ If a variable is not
$\Sigma$-essential in $t$, then it is called \emph{$\Sigma$-fictive}
 in $t$. $Fic(t,\Sigma)$ denotes the set of all $\Sigma$-fictive
variables  in $t.$

    The concept of
$\Sigma$-essential positions is a natural extension of
$\Sigma$-essential variables.

Let  $\mathcal{A}=\langle A,\mathcal{F}\rangle$ be
an algebra of type $\tau$, $t\in W_\tau(X_n)$, and let $p\in Pos(t)$.
 If $x_{n+1}\in
Ess(t(p;x_{n+1}),\mathcal{A})$, then the
  position $p\in Pos(t)$ is called
\emph{essential}  in $t$  with respect to   
$\mathcal{A}$. The set of all essential positions in  $t$  with
respect to  $\mathcal{A}$ is denoted by $PEss(t,\mathcal{A})$ (see Example \ref{ex1} below or Example 2.1 of \cite{sht1}).

When
a position $p\in Pos(t)$ is not essential in $t$ with respect to
$\mathcal{A}$, it is called \emph{fictive}  in $t$ with respect to
$\mathcal{A}$. The set of all fictive positions with respect to
$\mathcal{A}$ is denoted by $PFic(t,\mathcal{A}).$

If $x_{n+1}\in
Ess(t(p;x_{n+1}),\Sigma)$ the position $p\in Pos(t)$ is called
\emph{$\Sigma$-essential} in $t$ \cite{sht1}.  The set of $\Sigma$-essential
positions  in  $t$ is denoted by $PEss(t,\Sigma).$ When a position
is not $\Sigma$-essential in $t$ it is called \emph{
$\Sigma$-fictive}. $PFic(t,\Sigma)$ denotes the set of all
$\Sigma$-fictive positions  in $t.$

The set of $\Sigma$-essential subterms in $t$ is defined as follows:
$SEss(t,\Sigma)=\{r\in W_\tau(X) \mid  \Sigma\models r\approx sub_t(p),$ $\ p\in PEss(t,\Sigma)\}$.
 
So, a term  is a $\Sigma$-essential subterm of a term $t$, if it is $\Sigma$-equivalent to a subterm of $t$, whose root is located at 
$\Sigma$-essential positions in $t$.

Let $\Sigma$ be a set of identities of type $\tau$. Two terms $t$
and $s$ are called \emph{$\Sigma$-equivalent} (or \emph{$\Sigma$-equal}) if $\Sigma\models t\approx s$.
 
Let $t,r\in W_\tau(X)$ and let $\Sigma S_r^t=\{v\in Sub(t) \mid
\Sigma\models r\approx v \}$ be the set of all subterms of $t$ which
are $\Sigma$-equal to $r$.

Let $\Sigma P_r^t=\{p\in Pos(t) \mid sub_t(p)\in \Sigma S_r^t\}$ be
the set of all positions of subterms of $t$ which are $\Sigma$-equivalent
to $r$. Let $P_r^t=\{p_1,\ldots,p_m\}$ be the set of  all the
minimal elements in $\Sigma P_r^t$ with respect to the ordering
$\preceq$ in the set of positions, i.e. $p\in P_r^t$  if for each $q\in
\Sigma P_r^t$ we have $q\not\preceq p$.

\begin{definition}\label{d2.1}\cite{sht1}
\emph{ Term
$\Sigma$-composition} $t^\Sigma(r\leftarrow s)$ of $t$ and $r$ by $s$ is defined as follows

\begin{enumerate}
 \item[(i)]$t^\Sigma( r\leftarrow s)=t$\ \  if\  \ $P_r^t=\emptyset$;

\item[(ii)] $t^\Sigma(r\leftarrow s)= t(P_r^t;s)$\ \  if\ \ 
$P_r^t\neq\emptyset$.
\end{enumerate}
\end{definition}

\begin{lemma}\label{l2.1}
 If $\Sigma\models r\approx v$ then $t^\Sigma(r\leftarrow u)=t^\Sigma(v\leftarrow u)$.
\end{lemma}
\begin{proof}
 The lemma follows from the obvious equation  $P_r^t=P_v^t$ for each term $v\in W_\tau(X)$ with $\Sigma\models r\approx v$.
\end{proof}

\begin{example}\label{ex1}
  Let $\tau=(2)$ and let us consider the variety
$RB=Mod(\Sigma)$ of rectangular bands, where
\begin{equation*}~\label{eq2} \Sigma=\{x_1(x_2x_3) \approx (x_1x_2)x_3 \approx
x_1x_3,\ x_1x_1 \approx x_1\}.\end{equation*}
Let
$t=((x_1x_2)x_2)((x_1x_2)x_3)$, $r=x_1x_2$  and $s=x_4$.

 It is not difficult to see that the sets of $\Sigma-$essential
positions and subterms in $t$ are
\[PEss(t,\Sigma)=\{\varepsilon,\mbox{\small 1,11,111,2,22}\}\] and
\[SEss(t,\Sigma)=\{t,(x_1x_2)x_2,x_1x_2,x_1,(x_1x_2)x_3,x_3\}.\]

The $\Sigma-$essential and $\Sigma-$fictive positions in $t$ are
represented by large and small black circles, respectively in Figure
\ref{f-EssPos}. 
Next, we have 
\[\Sigma S_r^t=\{x_1x_2,(x_1x_2)x_2\},\quad \Sigma P_r^t=\{\mbox{\small 1,11,21}\}\quad and\quad P_r^t=\{\mbox{\small 1,21}\}.\]
Thus we have $t^\Sigma(r\leftarrow s)=x_4(x_4x_3)$ (see Figure \ref{fig2}).

\end{example}

\begin{figure}[ht]
\unitlength=1.00mm \special{em:linewidth 0.4pt}
\linethickness{0.4pt}
\begin{picture}(64.00,45.00)
\put(35.00,10.00){\circle*{2.00}}
\put(35.00,1.00){\makebox(0,0)[cc]{$t$}} \
\put(35.00,10.00){\line(-3,2){15.00}}
\put(35.00,10.00){\line(3,2){15.00}}
\put(20.00,20.00){\circle*{2.00}}
\put(50.00,20.00){\circle*{2.00}}

\put(50.00,20.00){\line(-1,1){10.00}}
\put(50.00,20.00){\line(1,1){10.00}}

\put(27.00,15.50){\line(1,2){4.00}}
\put(27.00,15.50){\circle*{2.00}}
\put(31.00,23.00){\circle*{1.00}}
\put(34.00,23.00){\makebox(0,0)[cc]{$x_2$}}
\put(32.00,20.00){\makebox(0,0)[cc]{\small 12}}
\put(24.00,15.50){\makebox(0,0)[cc]{\small 1}}

\put(20.00,20.00){\line(-1,2){5.00}}
\put(20.00,20.00){\line(1,2){5.00}}
\put(15.00,30.00){\circle*{2.00}}
\put(25.00,30.00){\circle*{1.00}}
\put(60.00,30.00){\circle*{2.00}}
\put(40.00,30.00){\circle*{1.00}}
\put(40.00,30.00){\line(-1,2){5.00}}
\put(40.00,30.00){\line(1,2){5.00}}
\put(35.00,40.00){\circle*{1.00}}
\put(45.00,40.00){\circle*{1.00}}
\put(17.00,20.00){\makebox(0,0)[cc]{\small 11}}
\put(37.00,30.00){\makebox(0,0)[cc]{\small 21}}
\put(64.00,29.00){\makebox(0,0)[cc]{$x_3$}}
\put(31.50,40.00){\makebox(0,0)[cc]{\small 211}}
\put(45.00,44.00){\makebox(0,0)[cc]{$x_2$}}
\put(25.00,34.00){\makebox(0,0)[cc]{$x_2$}}
\put(16.00,34.00){\makebox(0,0)[cc]{$x_1$}}
\put(11.00,30.00){\makebox(0,0)[cc]{\small 111}}
\put(21.00,30.00){\makebox(0,0)[cc]{\small 112}}
\put(36.00,45.00){\makebox(0,0)[cc]{$x_1$}}
\put(41.00,40.00){\makebox(0,0)[cc]{\small 212}}
\put(56.00,30.00){\makebox(0,0)[cc]{\small 22}}
\put(47.00,20.00){\makebox(0,0)[cc]{\small 2}}
\put(31.00,8.00){\makebox(0,0)[cc]{$\varepsilon$}}
\end{picture}

  \caption{$\Sigma-$essential positions in $t$.}
\label{f-EssPos}
\end{figure}
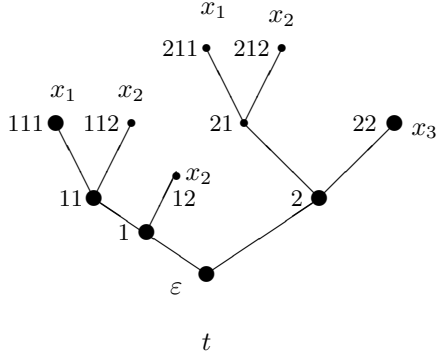

\begin{figure}[ht]
\unitlength=1.00mm \special{em:linewidth 0.4pt}
\linethickness{0.4pt}
\begin{picture}(64.00,45.00)
\put(35.00,10.00){\circle*{2.00}}
\put(35.00,1.00){\makebox(0,0)[cc]{$t^\Sigma(r\leftarrow s)$}} \
\put(35.00,10.00){\line(-3,2){15.00}}
\put(35.00,10.00){\line(3,2){15.00}}
\put(20.00,20.00){\circle*{2.00}}
\put(50.00,20.00){\circle*{2.00}}

\put(50.00,20.00){\line(-1,1){10.00}}
\put(50.00,20.00){\line(1,1){10.00}}

\put(60.00,30.00){\circle*{2.00}}
\put(40.00,30.00){\circle*{1.00}}

\put(17.00,20.00){\makebox(0,0)[cc]{\small 1}}
\put(24.00,20.00){\makebox(0,0)[cc]{$x_4$}}
\put(37.00,30.00){\makebox(0,0)[cc]{\small 21}}
\put(44.00,30.00){\makebox(0,0)[cc]{$x_4$}}
\put(64.00,29.00){\makebox(0,0)[cc]{$x_3$}}

\put(56.00,30.00){\makebox(0,0)[cc]{\small 22}}
\put(47.00,20.00){\makebox(0,0)[cc]{\small 2}}
\put(31.00,8.00){\makebox(0,0)[cc]{$\varepsilon$}}
\end{picture}

  \caption{$\Sigma-$composition of terms $t$ and $r$ by $s$.}
\label{fig2}
\end{figure}
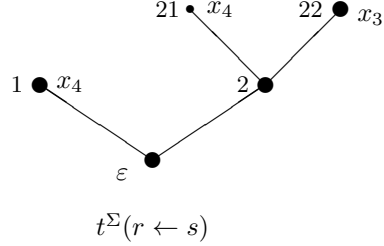

Terms are important tools in various areas, such as abstract data type
specifications, implementation of programming languages, automated deduction
etc. They can be used as
models for different structures in logic programming, term rewriting systems and
other computational procedures.

A term complexity measure or \emph{valuation} of a term is a  function
$Val:W_\tau(X)\to   \N$ if for some $c\in N$,  $Val(x_i)=c$ for
all $i\geq 1$ and $Val(t)\geq c$ for all $t\in W_\tau(X)$. The natural
number $c$ is called initial value of the valuation $Val$. It is often important for
applications that terms  be represented in forms with low complexity, including sometimes in normal forms.

Some common valuations are based on a linguistic point of view which counts
the number
of variables or the number of operation symbols occurring  in the term.

If $l_i$ denotes the number of occurrences of
the variable $x_i$ in the $n$-ary term $t$ then the  valuation $Len$
is called the \emph{length} of $t$ and it is defined as follows
$Len(t) = \sum \limits_{x_i\in var(t)} l_i .$  Its initial value is $1$.

The  \emph{depth} of a term $t$   is
defined  as follows: 
$Depth(x_i) = 0$ for $i=1,2,\ldots$  and
\\ \centerline{$Depth(f(t_1,\ldots,t_{n}))=max
\{Depth(t_1),\ldots,Depth(t_{n})\} +1.$}

 Let $t\in W_\tau(X)$ be a term and $Wv(t)=x_{i_1}\ldots x_{i_s}$ be the word of variables in $t$ which are written from the left to right, and let $st(t)=i_1\ldots i_s\in \N^*$ be the string of the indexes in  $Wv(t)$. The term $t$ is called {\it $\Sigma$-minimal} if for each $s\in W_\tau(X)$ with $\Sigma\models t\approx s$ it holds $Len(t)< Len(s)$ or $st(t)\prec_{lex} st(s)$ when $Len(t)=Len(s)$, where $\prec_{lex}$ is the lexicographical order in $\N^*$.
 
 Clearly, $\Sigma$-minimal terms are unique.

 For instance, let $t$ be the term defined in Example \ref{ex1}. Then we have 
  $Len(t)=6$, $Depth(t)=3$, $Wv(t)=x_1x_2x_2x_1x_2x_3$ and $st(t)=122123$. The $\Sigma$-minimal term corresponding to $t$ is $x_1x_3$.
 
 It is clear that $\Sigma$-minimal terms minimize the valuations $Depth$ and $Len$ in the sets of $\Sigma$-equal terms.

We need  some basic definitions from   universal algebra. More detailed
background about these notions can be found in the classic text \cite{bur}.

\begin{definition}
\cite{bur} A set $\Sigma$ of identities of type
$\tau$ is \emph{$D$-deductively closed} if it satisfies the following
axioms (some authors call them ``deductive rules", ``derivation
rules", ``productions", etc.):

\begin{enumerate}

\item[$D_1$] {(reflexivity)} $t\approx t\in \Sigma$ for
each term $t\in W_\tau(X)$;

\item[$D_2$] {(symmetry)}  $(t\approx s\in \Sigma)\  \Rightarrow\
s\approx t\in \Sigma$;

\item[$D_3$]  {(transitivity)} $(t\approx s\in \Sigma)\  \&\ (s\approx
r\in \Sigma)\  \Rightarrow\  t\approx r\in \Sigma$;

\item[$D_4$] {(term positional replacement)} \\
$(t\approx s\in \Sigma)\ \&\ (r\in W_\tau(X))\ \&\ (sub_r(p)=t)\
\Rightarrow\ r(p;s)\approx r\in\Sigma$;

\item[$D_5$] {(variable inductive substitution)} \\ $(t\approx s\in
\Sigma)\ \&\ (r\in W_\tau(X))\   \Rightarrow\ t(x\leftarrow
r)\approx s(x\leftarrow r)\in\Sigma$.

\end{enumerate}
\end{definition}
For any set $\Sigma$ of identities,  the smallest $D$-deductively
closed set containing $\Sigma$ is called the $D$-closure of $\Sigma$
and it is denoted by $D(\Sigma).$

The first three deductive rules make $D(\Sigma)$ into an equivalence relation, the fourth makes it a congruence, and the last rule says $D(\Sigma)$ is a fully invariant congruence.

Let $\Sigma$ be a set of identities of type $\tau.$ For $t\approx
s\in Id(\tau)$ we say $\Sigma$ proves  $t\approx s$ and write $\Sigma\vdash t\approx s$ if there is a sequence of
identities
($D$-deduction) $t_1\approx s_1,\ldots,t_n\approx s_n$, such that
each identity belongs to $\Sigma$ or is a result of applying any of
the derivation rules $D_1-D_5$  to previous identities in the
sequence and the last identity $t_n\approx s_n$ is $t\approx s.$

It is well known that $\Sigma\vdash t\approx s\ \iff \ \Sigma\models t\approx
s$.

In \cite{sht1} a variation of the derivation rules  $D_1 - D_5$ is given, which is used to define a globally
invariant congruence. 
\begin{definition}\label{d2} \cite{sht1}
 A set   $\Sigma$ of identities  is \emph{$\Sigma R$-deductively closed} if
it satisfies the rules $D_1,D_2,D_3,D_5$ and

\begin{enumerate}

\item[$\Sigma R_1$] \emph{($\Sigma$  replacement)}
\[{\left(\begin{array}{c}
r,t,s,u\in W_\tau(X)\ \&\ (t\approx s\in\Sigma)\ \&\   \\
                 r\in SEss(t,\Sigma)\cap SEss(s,\Sigma)
\end{array}\right)}\Rightarrow t^\Sigma(r\leftarrow u)\approx
s^\Sigma(r\leftarrow u)\in\Sigma.\]
\end{enumerate}
\end{definition} 

For any  set $\Sigma$ of identities, the smallest $\Sigma
R$-deductively closed set containing $\Sigma$ is called  the \emph{$\Sigma
R$-closure} of $\Sigma$  and it is denoted by $\Sigma R(\Sigma).$

  $\Sigma R$ is a closure operator which implies that:
\begin{enumerate}  
 
\item[(1)] $\Sigma R(\Sigma R(\Sigma))= \Sigma R(\Sigma)$ and 

\item[(2)] for each $\Delta\subseteq \Sigma$,  if $\Delta\vdash t\approx s$ then $t\approx s\in \Sigma R(\Sigma)$. 
\end{enumerate}

A set $\Sigma\subseteq Id(\tau)$ is called \emph{a globally invariant congruence} if it is $\Sigma
R$-deductively closed.  In \cite{sht1} it is proved that each globally invariant congruence is a fully invariant congruence.

A variety $V$ of type $\tau$ is called \emph{ stable} if $\Sigma =Id(V)$ is $\Sigma
R$-deductively closed.

\section{Stable varieties of semigroups}\label{sec4}

  We are going to describe all stable varieties of semigroups in an analogy to the  solid varieties \cite{lpol2,lpol}, using  some fundamental results in semigroup theory,  as presented in  \cite{evans,hov}. 

Let us  agree,  throughout the rest of the paper to write $f(x_1,x_2)$ as $(x_1x_2)$ or  $x_1x_2$. 

The following identities of type $(2)$  are important for the achievement of our aim:
\begin{equation}\label{eq.1}  x_1x_2x_3 \approx
x_ix_j\quad for\quad  i,j\in\{1,2,3\}. \end{equation}
They allow us to define a special class of varieties of semigroups. Let $i$ and $j$  be two natural numbers from the set $\{1,2,3\}$. Then we consider the following variety of semigroups: 
\[V_{ij}=Mod(\{(x_1x_2)x_3\approx x_1(x_2x_3),\  x_1x_2x_3\approx x_{i}x_{j}\}).\]

Let $t\in W_\tau(X)$ be a term and $Wv(t)=x_{i_1}\ldots x_{i_s}$ be the string of the variables in $t$. We use the notations $first(t)=x_{i_1}, second(t)=x_{i_2},\ldots,$ $s$-$th(t)=x_{i_s}$. Often, the last variable $x_{i_s}$ is denoted by $last(t)$ or $rightmost(t)$. Also, some authors write $first(t)$ as $leftmost(t)$.  For instance, if   $t=x_3x_1x_2x_2$ then we have $first(t)=leftmost(t)=x_3$, $second(t)=x_1$, $third(t)=x_2$ and $fourth(t)=last(t)=rightmost(t)=x_2$. 

\begin{lemma}\label{l3.1}
The varieties $LZ$ (of Left-Zero-semigroups), $RZ$ (of Right-Zero-semigroups), $Z$ (of Zero-semigroups), and the varieties  $V_{ij}$
 with $1\leq i< j\leq 3$ are stable.
\end{lemma}

\begin{proof} Let $\mathcal V\in \{LZ, RZ, Z, V_{12}, V_{13}, V_{23}\}$ and $\Sigma=Id(\mathcal V)$.
Since $Id(\mathcal V)$ is a fully invariant congruence, it satisfies the derivation  rules $D_1$, $D_2$, $D_3$, $D_4$ and $D_5$. We have to prove that $\Sigma R_1$ is also satisfied in $\mathcal V$, i.e.

\begin{equation}\label{eq.2}
\Sigma\models t^\Sigma(r\leftarrow u)\approx s^\Sigma(r\leftarrow
u),
\end{equation}
when $\Sigma\models t\approx s$,  $r\in SEss(t,\Sigma)\cap SEss(s,\Sigma)$ and $u\in W_\tau(X)$.

 Let
$t,s,r,u\in W_\tau(X)$ be four   terms such that  $\Sigma\models
t\approx s,$ 
$r\in SEss(t,\Sigma)\cap SEss(s,\Sigma)$ and $u\in W_\tau(X)$.

To prove (\ref{eq.2}), let us observe the following two common facts:

First, if $P_r^t=\{\varepsilon\}$ or $P_r^s=\{\varepsilon\}$ then $r\in SEss(t,\Sigma)\cap SEss(s,\Sigma)$ implies $\Sigma\models t\approx r$ and $\Sigma\models s\approx r$. Thus we have \[t^\Sigma(r\leftarrow u)=s^\Sigma(r\leftarrow u)=u,\]
which shows that (\ref{eq.2}) is satisfied.

Second,  if $Depth(t)=0$ then we have $t=x_i$ for some variable $x_i\in X$  and hence $r\in SEss(t,\Sigma)\cap SEss(s,\Sigma)$ implies $\Sigma\models r\approx x_i$.
Now,  (\ref{eq.2}) is satisfied, according to $D_5$.
 
Next, we assume that  $\Sigma\models
t\approx s,$
$r\in SEss(t,\Sigma)\cap SEss(s,\Sigma)$, $1\leq Depth(t)\leq Depth(s)$ and $\Sigma\not\models
t\approx r$.
\vspace{.5cm}

\noindent
{\bf Claim 1.} $LZ$ and $RZ$ are stable varieties.
\vspace{.5cm}

\noindent
Let us consider the variety $\mathcal V=LZ$. Then $\Sigma \models x_1x_2\approx x_{1}$ and  clearly $\Sigma\models w\approx first(w)$ for all terms $w$. Consequently, $\Sigma\models
t\approx s$ and 
$r\in SEss(t,\Sigma)\cap SEss(s,\Sigma)$ imply
\[first(t)=first(s)=first(r).\]
Thus we have $\Sigma\models t\approx s\approx r$ and hence  (\ref{eq.2}) is satisfied,
which shows that $LZ$ is stable.

The variety $RZ$ is stable by dual arguments.
\vspace{.5cm}

\noindent
{\bf Claim 2.} $Z$ is a stable variety.
\vspace{.5cm}

\noindent
We have $\Sigma \models x_1x_2\approx x_{3}x_4$ and  clearly $P_r^t=P_r^s=\{\varepsilon\}$ for all terms $t,s$ and $r$ with $\Sigma\models
t\approx s$ and 
$r\in SEss(t,\Sigma)\cap SEss(s,\Sigma)$,   which proves the stability of $Z$.
\vspace{.5cm}

\noindent
{\bf Claim 3.} $V_{12}, V_{13}$ and $V_{23}$ are stable varieties.
\vspace{.5cm}

\noindent
We shall show that the variety  $\mathcal V=V_{12}$ is stable. Hence, we have  \[\Sigma\models (x_1x_2)x_3\approx x_1(x_2x_3) \quad\mbox{and}\quad \Sigma\models x_1x_2x_3\approx x_{1}x_{2}.\]
Since $Depth(t)\geq 1$, we have $\Sigma\models t\approx first(t)second(t)$. Let us assume, with no loss of generality  that $x_1=first(t)$ and $x_2=second(t)$.

Next, $\Sigma\not\models t\approx r$ and $r\in SEss(t,\Sigma)\cap SEss(s,\Sigma)$ implies $\Sigma\models r\approx x_1$ or $\Sigma\models r\approx x_2$. Without loss of generality let us assume that $\Sigma\models r\approx x_1$. According to Lemma \ref{l2.1}, we have
\[\Sigma\models t^\Sigma(r\leftarrow u)\approx t^\Sigma(x_1\leftarrow
u).\]
Since $x_1$ is a variable it is easy to see that 
\[\Sigma\models t^\Sigma(x_1\leftarrow u)\approx t(x_1\leftarrow
u).\]
Hence for satisfaction of (\ref{eq.2}) we need  
\[\Sigma\models t(x_1\leftarrow u)\approx s(x_1\leftarrow
u),\]
which follows from  $D_5$. Consequently, $V_{12}$ is a stable variety. The proof that $V_{13}$ and $V_{23}$ are stable varieties is left to the reader.
\end{proof}
\begin{remark}\label{r2}
 Let us consider the variety $V_{21}=Mod(\Sigma)$, where  \[\Sigma = \{(x_1x_2)x_3\approx x_1(x_2x_3),\  x_1x_2x_3\approx x_{2}x_{1}\}.\]
Then  we have 
\begin{eqnarray*}\label{eq.7} 
 \Sigma\models x_1x_2\approx (x_2x_1)x_3\approx ((x_1x_2)x_4)x_3\approx x_1(x_2x_4)x_3\approx x_1((x_4x_2)x_5)x_3\approx \nonumber\\
 \approx ((x_1x_4)x_2)x_5x_3\approx x_4((x_1x_5)x_3)\approx x_4(x_5x_1)\approx ((x_4x_5)x_1)\approx x_5x_4.
\end{eqnarray*}
 Hence $V_{21}=Z$. 
Using similar or dual arguments one can show that  $V_{31}=V_{32}=Z$.
\end{remark}

\begin{proposition}\label{p1}
The varieties  of semigroups $V_{ii}$,
for $i\in\{1,2,3\}$  are not stable. 
\end{proposition}
\begin{proof}
 We shall prove that $V_{11}= Mod(\{(x_1x_2)x_3\approx x_1(x_2x_3),\  x_1x_2x_3\approx x_{1}x_{1}\})$ is not a stable variety.
 
 Let us consider the following terms $t=(x_1x_2)x_3$, $s=(x_1x_2)x_4$ and $r=x_1x_2$. Clearly, $\Sigma\models t\approx s$. Then we have $t(${\small 1}$;x_1)=x_1x_3$ and $t(${\small 1}$;x_3)=x_3x_3$ and $\Sigma\not\models x_3x_3\approx x_1x_3$. Since $\Sigma\not\models x_1x_2\approx (x_1x_2)x_3$ we have {\small 1}$\in P_r^t$. In an analogous way we obtain {\small 1}$\in P_r^s$. Next, we have $t^\Sigma(r\leftarrow x_3)=x_3x_3$ and $s^\Sigma(r\leftarrow x_3)=x_3x_4$ which shows that $V_{11}$ is not stable. 
 
 In a similar way one can prove that $V_{22}$  and   $V_{33}$ are  not stable varieties.
\end{proof}

\begin{lemma}\label{l3}
 The varieties of semigroups
\[V_1=Mod(\{(x_1x_2)x_3\approx x_1(x_2x_3),\  x_1x_2\approx x_1x_3\})\] and
  \[V_3=Mod(\{(x_1x_2)x_3\approx x_1(x_2x_3),\   x_1x_3\approx x_2x_3\})\]
 are stable.
\end{lemma}
\begin{proof} We shall prove that $V_1$ is stable.
 Let $t,s$ and $r$ be terms for which $\Sigma\models t\approx s$ and $r\in SEss(t,\Sigma)\cap SEss(s,\Sigma)$.
 
 If $Depth(t)=0$ then $Depth(s)=0$ and (\ref{eq.2}) is clearly satisfied.
 
 If $Depth(t)=1$ and $Depth(s)=1$ then $first(t)=first(s)$ and $r=first(t)$ or $\Sigma\models t\approx r$, and (\ref{eq.2}) is obvious, again.

 If $Depth(t)\geq 2$ then $\Sigma\models t\approx s\approx first(t)x_2$ for an arbitrary  variable $x_2\in X$. This implies that $\Sigma\models r\approx first(t)$ or $\Sigma\models t\approx r$. In both cases (\ref{eq.2}) is satisfied.

 By dual arguments it follows that $V_3$ is stable.
 
\end{proof}

\begin{remark}\label{r1} 
Since $\{x_1x_2\approx x_1x_3\}\models x_1x_2x_3\approx x_{1}x_{1}$ it follows that $V_{1}\subseteq V_{11}$ and  by dual arguments we have $V_{3}\subseteq V_{33}$.
\end{remark}

\begin{lemma}\label{l3.2}
Let $\mathcal V=Mod(\Sigma)$ be a stable variety of semigroups. If \    
 $\Sigma\models x_1x_2x_3\approx x_1x_2x_4$  then 
 $\Sigma$   proves at least one  identity among (\ref{eq.1}) with $1\leq i< j\leq 3$.
 
\end{lemma}
 
\begin{proof}  Let us consider the following terms $t=(x_1x_2)x_3$, $s=(x_1x_2)x_4$, $r=x_1x_2$ and $u=x_1$. Clearly, $\Sigma\models t\approx s$.
 
If $r\notin SEss(t,\Sigma) \cap  SEss(s,\Sigma)$  then we are done because
$\Sigma\models (x_1x_2)x_3\approx x_5x_3$ and from $D_5$ we have $\Sigma\models (x_1x_2)x_3 \approx x_1x_3$. 

Let $r\in SEss(t,\Sigma)\cap SEss(s,\Sigma)$. If  $\Sigma\models t\approx r$ then
$\Sigma\models s\approx r$ and we are done again,  because of $\Sigma\models (x_1x_2)x_3\approx x_1x_2$.

Next, assume that  $\Sigma \not\models t\approx r$. Then $\Sigma \not\models s\approx
r$ and  $P_r^t=\{${\small 1}$\}$ and $P_r^s=\{${\small 1}$\}$. From   
$\Sigma \models t\approx s$  and $\Sigma R_1$ we obtain 
\[\Sigma\models t^\Sigma(r\leftarrow x_1)\approx  s^\Sigma(r\leftarrow x_1)\quad\mbox{and}\quad\Sigma \models x_1x_3\approx x_1x_4.\]

According to $D_5$ we can replace  $x_4$
by $x_2x_3$ in the last identity, and hence $\Sigma \models
 x_1x_3\approx x_1(x_2x_3)$.
 \end{proof}

In a similar way, one can show that if $\Sigma \models x_1x_2x_3 \approx
x_4x_2x_3$ or $\Sigma \models x_1x_2x_3 \approx
x_1x_4x_3$ then  $\Sigma$   proves at least one  identity among (\ref{eq.1}) with $1\leq i<j\leq 3$.

\begin{lemma}\label{l3.3}
If $\mathcal V=Mod(\Sigma)$ is  a stable  variety of semigroups 
 then  
$ \Sigma \models x_1x_1x_1\approx
x_1x_1.$

\end{lemma}
\begin{proof}
If $\Sigma$  proves at least one  identity among (\ref{eq.1}) with  $1\leq i<j\leq 3$ then  $\Sigma \models x_1x_1x_1\approx
x_1x_1$ is clear.

Let us assume that $\Sigma$ does not prove any  identity among (\ref{eq.1}) with  $1\leq i<j\leq 3$. 
 We shall prove the lemma by  considering cases:
\vspace{.5cm}

\noindent
{\bf Case A.}  $\Sigma\not\models (x_1x_2)(x_1x_2)\approx x_1x_2$. 
\vspace{.5cm}

\noindent 
Let us put $t=((x_1x_2)x_1)x_2$, $s=(x_1x_2)(x_1x_2)$,
$r=x_1x_2$ and $u=x_3$. Clearly $\Sigma\models t\approx s$. If we suppose
that $r\notin SEss(t,\Sigma)$
then $\Sigma\models (x_3x_1)x_2\approx (x_4x_1)x_2$ which
contradicts Lemma \ref{l3.2} and hence $r\in  SEss(t,\Sigma)$. If we suppose
that $r\notin SEss(s,\Sigma)$ then $\Sigma\models x_3x_1\approx
x_3x_4\approx x_3(x_1x_2)$  which is a contradiction. Hence  $r\in 
SEss(s,\Sigma)$.

Then we have
$t^\Sigma(r\leftarrow u)=(x_3x_1)x_2$ and $s^\Sigma(r\leftarrow
u)=x_3x_3$. Hence $\Sigma\models (x_3x_1)x_2\approx x_3x_3$ and after replacing $x_3$ and $x_2$ by $x_1$,  we obtain 
$\Sigma\models x_1x_1x_1\approx x_1x_1$.
\vspace{.5cm}

\noindent
{\bf Case B.} $\Sigma\models (x_1x_2)(x_1x_2)\approx x_1x_2$. 
\vspace{.5cm}

\noindent
The associative law and  $D_5$ imply 
 \begin{equation}\label{eq.3}\Sigma\models ((x_1x_1)x_1)x_1\approx
(x_1x_1)(x_1x_1)\approx x_1x_1.\end{equation}
Let us put $t=((x_1x_2)x_3)(x_1x_2)$,
$s=(((x_1x_2)x_3)x_1)x_2$,
$r=x_1x_2$ and $u=x_4$.  Clearly $\Sigma\models t\approx s$. 

If  $\Sigma\models t\approx r$ then we are done  after replacing
$x_2$ and $x_3$ by $x_1$ in $\Sigma\models t\approx s$.

If we suppose that $r\notin SEss(t,\Sigma)$ then  $\Sigma\models
(x_4x_3)x_5\approx (x_6x_3)x_7$ which contradicts Lemma \ref{l3.2}.
Hence $r\in SEss(t,\Sigma)$. 

Let us assume that $r\in SEss(t,\Sigma)\cap SEss(s,\Sigma)$ then from  (\ref{eq.2}) we obtain
$\Sigma\models
(x_4x_3)x_4\approx ((x_4x_3)x_1)x_2.$
Replacing $x_2$, $x_3$ and $x_4$ by $x_1$ and using
(\ref{eq.3}), we obtain   $\Sigma\models x_1x_1\approx
(x_1x_1)x_1$.

Assume that $r\in SEss(t,\Sigma)\setminus SEss(s,\Sigma)$
then   
\[\Sigma\models (((x_1x_2)x_3)x_1)x_2\approx
((x_4x_3)x_1)x_2.\]

Again,  replacing $x_2$, $x_3$ and $x_4$ by $x_1$ and  using (\ref{eq.3}),
we obtain   $\Sigma\models x_1x_1\approx
x_1x_1x_1$.
\end{proof}
\begin{theorem}\label{t3.1}
If a variety $\mathcal V$ of semigroups  is stable then  $\mathcal V\subseteq V_{12}$ or $\mathcal V\subseteq V_{13}$, or $\mathcal V\subseteq V_{23}$.
\end{theorem}
\begin{proof}
 Let $\mathcal V=Mod(\Sigma)$ be a stable variety of semigroups. 

First,  let
$
 \Sigma \not\models ((x_1x_2)x_2)x_3 \approx
((x_1x_2)x_2)x_4
$
and let us put $t= ((x_1x_2)x_2)x_3$, $s=
(((x_1x_2)x_2)x_2)x_3$,  $r=(x_1x_2)x_2$ and $u=x_4$. Lemma
\ref{l3.3} implies that $\Sigma\models t\approx s$.

If  $r\notin SEss(t,\Sigma)$ then   $\Sigma\models
x_1x_3\approx x_2x_3$ and according to $D_5$ we have
$\Sigma\models x_1x_3\approx (x_1x_2)x_3$. 

If
$r\notin SEss(s,\Sigma)$ then  $\Sigma\models (x_1x_2)x_3\approx
(x_4x_2)x_3$ and hence we are done because of Lemma \ref{l3.2}. 

If we suppose that $\Sigma\models r\approx t$, i.e.  $\Sigma\models (x_1x_2)x_2
\approx
((x_1x_2)x_2)x_3$, then $D_3$ implies $\Sigma\models (x_1x_2)x_2
\approx
((x_1x_2)x_2)x_4$ which contradicts our assumption.
Hence $\Sigma\not\models r\approx t$, $P_r^t=\{\mbox{\small 1}\}$, $P_r^s=\{\mbox{\small 11}\}$ and  (\ref{eq.2})
implies $\Sigma \models x_4x_2\approx (x_4x_2)x_3$. 

Second, assume that
$
 \Sigma\models ((x_1x_2)x_2)x_3 \approx
((x_1x_2)x_2)x_4$.
Lemma \ref{l3.3} implies  $\Sigma\models ((x_1x_2)x_2)x_3 \approx
((x_1x_2)x_2)x_2\approx (x_1x_2)x_2$.
Let us consider the following terms $t= (x_1x_2)x_2$, $s= ((x_1x_2)x_2)x_3$,
$r=x_1x_2$ and $u=x_4$.
Clearly $\Sigma\models t\approx s$.

If  $r\notin SEss(t,\Sigma)$ then   $\Sigma\models
x_1x_2\approx x_3x_2$ and hence $\Sigma\models x_1x_2\approx
(x_1x_3)x_2$.

If $r\notin SEss(s,\Sigma)$ then  $\Sigma\models (x_1x_2)x_3\approx
(x_4x_2)x_3$ and hence we are done because of Lemma \ref{l3.2}.

 Let
$r\in SEss(t,\Sigma)\cap SEss(s,\Sigma)$.

If  $\Sigma\models t\approx r$ then $\Sigma\models s\approx r$, i.e.
$\Sigma\models ((x_1x_2)x_2)x_3\approx x_1x_2\approx
(x_1x_2)x_3$. 

If $\Sigma\not\models t\approx
r$ then $P_r^t=\{\mbox{\small 1}\}$ and $P_r^s=\{\mbox{\small 11}\}$. Thus from 
(\ref{eq.2}) we obtain $\Sigma \models x_4x_2\approx (x_4x_2)x_3$, which
  completes the proof.
\end{proof}
\begin{theorem}\label{t3.2}
Let $\mathcal V$ be a variety of semigroups. Then $\mathcal V$ is stable if and only if $\mathcal V$
is one of the following ten varieties:
\begin{enumerate}
 \item[] $\mathcal V_1=Mod(\{(x_1x_2)x_3\approx x_1(x_2x_3),\  x_1x_2x_3\approx x_1x_3\})$;
 \item[] $\mathcal V_2=Mod(\{(x_1x_2)x_3\approx x_1(x_2x_3),\  x_1x_2x_3\approx x_1x_2\})$;
 \item[] $\mathcal V_3=Mod(\{(x_1x_2)x_3\approx x_1(x_2x_3),\  x_1x_2x_3\approx x_2x_3\})$;
   \item[] $\mathcal V_4=Mod(\{(x_1x_2)x_3\approx x_1(x_2x_3),\  x_1x_2\approx x_1x_3\})$;
 \item[] $\mathcal V_5=Mod(\{(x_1x_2)x_3\approx x_1(x_2x_3),\  x_1x_2\approx x_3x_2\})$;
 \item[] $\mathcal V_{6}=RB$ the variety of rectangular bands;
 \item[] $\mathcal V_{7}=LZ$ the variety of Left-Zero-semigroups;
 \item[] $\mathcal V_{8}=RZ$ the variety of Right-Zero-semigroups;
 \item[] $\mathcal V_{9}=Z$ the variety of Zero-semigroups;
 \item[] $\mathcal V_{10}=TR$ the trivial variety.
\end{enumerate}
All these varieties are pairwise distinct.
\end{theorem}
\begin{proof}
First, assume that $\mathcal V$ is a stable variety of semigroups. Then from Theorem \ref{t3.1} we have $\mathcal V=Mod(\Sigma)$ for some set of identities  $\Sigma$, which  proves at least one identity among (\ref{eq.1}) with $1\leq i< j\leq 3$. 
Then 
\[\mathcal V\in\{\mathcal V_{1},\mathcal V_{2},\mathcal V_{3},\mathcal V_{4},\mathcal V_{5},\mathcal V_{6},\mathcal V_{7},\mathcal V_{8},\mathcal V_{9},\mathcal V_{10}\}\]
follows by the following well known facts (see \cite{evans}):

{\it Fact 1.} The non-trivial proper subvarieties of $\mathcal V_2$ are $\mathcal V_4$, $LZ$ and $Z$.

{\it Fact 2.} The non-trivial proper subvarieties of $\mathcal V_3$ are $\mathcal V_5$, $RZ$ and $Z$.

{\it Fact 3.} The non-trivial proper subvarieties of $\mathcal V_1$ are $RB$, $\mathcal V_4$, $\mathcal V_5$, $LZ$,  $RZ$ and $Z$.

Second, the varieties $\mathcal V_1,\mathcal V_2, \mathcal V_3, \mathcal V_7, \mathcal V_8$ and $\mathcal V_9$ are stable according to Lemma \ref{l3.1}, the varieties  $\mathcal V_4$ and $\mathcal V_5$ are stable according to Lemma \ref{l3}. The stability of  $\mathcal V_6$  is proved in \cite{sht1}. The variety  $\mathcal V_{10}$ is obviously stable.

 Finally, it is a well known fact that all these varieties are pairwise distinct \cite{evans}.
 
All stable varieties of semigroups are shown in Figure \ref{fig3}.
\end{proof}
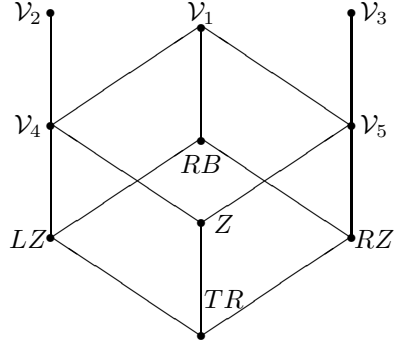
\begin{figure}[ht]
\centering

\unitlength 1mm 
\linethickness{0.4pt}
\begin{picture}(110.006,45.001)(0,0)

\put(40,0){\line(3,2){20}}
\put(40,0){\line(-3,2){20}}

\put(40,0){\circle*{1}}
\put(40,0){\line(0,1){15}}
\put(40,15){\circle*{1}}

\put(40,15){\line(3,2){20}}
\put(40,15){\line(-3,2){20}}

\put(20,13){\line(0,1){15}}
\put(60,13){\line(0,1){15}}

\put(20,13){\circle*{1}}
\put(60,13){\circle*{1}}

\put(20,28){\circle*{1}}
\put(60,28){\circle*{1}}

\put(20,28){\line(0,1){15}}
\put(60,28){\line(0,1){15}}

\put(20,43){\circle*{1}}
\put(60,43){\circle*{1}}
\put(60,28){\line(-3,2){20}}
\put(20,28){\line(3,2){20}}

\put(60,13){\line(-3,2){20}}
\put(20,13){\line(3,2){20}}

\put(40,26){\circle*{1}}
\put(40,26){\line(0,1){15}}
\put(40,41){\circle*{1}}

\put(40, 43){\makebox(0,0)[cc]{$\mathcal V_1$}}
\put(40, 23){\makebox(0,0)[cc]{$RB$}}

\put(17, 43){\makebox(0,0)[cc]{$\mathcal V_2$}}

\put(63,43){\makebox(0,0)[cc]{$\mathcal V_3$}}

\put(17, 28){\makebox(0,0)[cc]{$\mathcal V_4$}}
\put(17, 13){\makebox(0,0)[cc]{$LZ$}}
\put(63, 28){\makebox(0,0)[cc]{$\mathcal V_5$}}
\put(63, 13){\makebox(0,0)[cc]{$RZ$}}

\put(43, 15){\makebox(0,0)[cc]{${Z}$}}
\put(43, 5){\makebox(0,0)[cc]{$TR$}}

\end{picture}
\caption{Stable varieties of semigroups.}
\label{fig3}

\end{figure}

\section{Abstract reduction systems and deduction of identities}\label{sec3}

A Term Rewriting System (TRS) for deductions on identities is
a pair $(\tau, \mathcal R)$ of a type and a set of reduction (rewrite) rules,
which are binary relations on $W_\tau(X)$ written as  $t\rightarrow r$.

Our aim is to use TRS and apply their well developed tools to investigate the stability of several varieties of groupoids.   
For this purpose we  consider TRS as  Abstract Reduction Systems (ARS).

An ARS is a structure $\mathcal W=\langle W_\tau(X), (\rightarrow_i)_{i\in I},\Sigma
\rangle$, where $(\rightarrow_i)_{i\in I}$ is a family of binary relations on
$W_\tau(X)$, called \emph{reductions or rewrite relations}. 
For a reduction $\rightarrow_i$ the transitive and reflexive closure is denoted
by $\twoheadrightarrow_i$.
A term $r\in W_\tau(X)$ is a \emph{normal form} if there is no $v\in W_\tau(X)$
such that $r\rightarrow_i v$.

TRS, and in particular ARS, play an important role in various areas such as abstract data
type specification, functional programming, automated deductions, etc.  For more
detailed information about TRS we refer to J. W. Klop and Roel de Vrijer
\cite{klo}. The concepts and properties of ARS also apply to other rewrite
systems such  as string rewrite systems (Thue systems),  tree rewrite systems, graph
grammars, etc.

Many computations,
constructions, processes, translations, mappings and so on, can be modeled
as stepwise transformations of objects known as rewriting systems.  In all  
different branches of rewriting two basic concepts occur,
 known as termination (guaranteeing the
existence
of normal forms) and confluence (securing the uniqueness of normal
forms).

 Let us consider the  ARS  $\mathcal W=\langle W_{(2)}(X),
\{\rightarrow_{R}\},\Sigma\rangle$ determined by the following reduction:

\[t\rightarrow_{R} r\ \stackrel{def}{\iff}\ r=t(p;u)\] 
where $s=sub_t(p)$, $\Sigma\models s\approx u$ and $u$ is $\Sigma$-minimal.

According to $D_4$ we have that if $t\rightarrow_R r$ then $\Sigma\models t\approx r$.

Our intention is to reduce the terms in an identity to  normal forms and then
implement the deductive rules on these normal forms, preferably with low complexity terms.

First, we are interested in existence and uniqueness of normal forms for the
reduction  $\rightarrow_{R}$. 

A reduction $\rightarrow$ has the \emph{unique normal form property} (UN) if whenever
$t,r\in W_\tau(X)$ are normal forms and $\Sigma\models t\approx r$ then $t=r$.

We are going to prove that $\rightarrow_{R}$ is UN when $\Sigma$ determines the variety of idempotent groupoids or consists of identities as from (\ref{eq.1}). This
we shall do  using   Newman's Lemma (Theorem 1.2.1. \cite{klo}).

A reduction $\rightarrow$ is \emph{terminating} (or \emph{strongly normalizing} SN) if every reduction sequence
$t\rightarrow t_1\rightarrow t_2\ldots$ eventually must terminate.

A reduction $\rightarrow$ is \emph{weakly
confluent} (or   \emph{has weakly Church-Rosser property} WCR) if $t\rightarrow r$ and $t\rightarrow v$ imply
that there is  $w\in W_\tau(X)$ such that $r\twoheadrightarrow w$ and $v\twoheadrightarrow
w$.

\begin{theorem}\label{t4.1}
The reduction $\rightarrow_{R}$ is  SN  and WCR.
\end{theorem}
\begin{proof}
(SN) Clearly, if $t \rightarrow_{R} r$ then  $Len(t)\geq Len(r)$ or $st(t)\prec_{lex} st(r)$ when $Len(t)=Len(r)$.

Since the lengths $Len(z)$  of the terms $z$ in any
reduction sequence   decrease and strings $st(z)$ strongly decrease it follows that the sequence eventually must
terminate, i.e.  the reduction is terminating.

(WCR) Let $t$ be a term, $p,q\in Pos(t)$, $s=sub_t(p)$, $w=sub_t(q)$, $r=t(p;u)$ and $v=t(q;z)$ where $u$ and $z$ are $\Sigma$-minimal. 

If $p\prec q$ then we have
$t\rightarrow_R r\rightarrow_R w$. If  $q\prec p$ then $t\rightarrow_R w\rightarrow_R r$  which shows that 
reduction $\rightarrow_R$ is WCR in these two cases.

Let  $p\not\prec q$ and  $q\not\prec p$ and 
let $y$ be the $\Sigma$-minimal term with $\Sigma\models t\approx y$. Then we have 
\[t\rightarrow_R w\rightarrow_R r\twoheadrightarrow_R y\]
and \[t\rightarrow_R r\rightarrow_R w\twoheadrightarrow_R y.\]

\end{proof}

\begin{corollary}\label{c4.1}
The reduction $\rightarrow_{R}$ is UN.
\end{corollary}

\begin{proof}
 This follows from Newman's Lemma, which states that  WCR $\&$ SN $\Rightarrow$ UN (see
Theorem 1.2.1. \cite{klo}).  
\end{proof}

For each term $t\in W_\tau(X)$ we denote by $Red(t)$ the normal form
obtained from $t$ under the reduction $\rightarrow_{R}$.

\begin{corollary}\label{l4.1}
$\Sigma\models t\approx Red(t)$  for any term $t\in W_\tau(X)$. 
\end{corollary}

It is easy to see that the normal form operator $Red$ minimizes the valuations $Len$ and $Depth$.

\section{Stable varieties of groupoids}\label{sec5}

We are going to  study 
stable varieties of groupoids.
Let us note that if 
 $\Sigma=\emptyset$ then $Mod(\Sigma)$ is a stable variety.

 First, we consider the variety of idempotent groupoids. 

Note that if $t\in W_\tau(X)$ and $s\in Sub(Red(t))$ then there is $r\in Sub(t)$ such that
$\Sigma\models r\approx s$ and if $t=t_1t_2$ then $\Sigma\models Red(t)\approx
Red(t_1)Red(t_2)$.

\begin{lemma}\label{l5.2}
 If
$\Sigma=\{x_1x_1  \approx
x_1\}$  then  \begin{equation}\label{eq.4} \Sigma\models Red(t^\Sigma(r\leftarrow u))\approx
Red(t)^\Sigma(r\leftarrow u)\end{equation}
 for every $r,t,u\in W_\tau(X)$. 
\end{lemma}
\begin{proof} 
 Let
$t,r,u\in W_\tau(X)$ be three  terms. We shall proceed by induction on $Depth(t)$.

If $Depth(t)=0$ then $t=x_i$ for some natural number $i$. Then $Red(t)=x_i$ and
it is obvious that   (\ref{eq.4}) is satisfied.

Let us assume that for some natural number $k\geq 2$, if $Depth(t)<k$ then (\ref{eq.4}) is satisfied for $t$.

Let $Depth(t)=k$  and $t=t_1t_2$.
Let $r\in Sub(t)$. 
If $\Sigma\models r\approx t$ then $\Sigma\models r\approx Red(t)$ and  we have  $\Sigma\models r \approx Red(t).$
Hence $Red(t^\Sigma(r\leftarrow u))=u$ and $Red(t)^\Sigma(r\leftarrow u)=u$, which
proves (\ref{eq.4}).

Let $r\in Sub(t)$ and $\Sigma\not\models r\approx t$. 
If $\Sigma S_r^t=\emptyset$ then clearly $\Sigma S_{r}^{Red(t)}=\emptyset$ and (\ref{eq.4}) is
obviously satisfied in this case.

Assume that $\Sigma S_r^t\neq\emptyset$. By the inductive assumption we
have 
\[\Sigma\models Red((t_i^\Sigma(r\leftarrow u)\approx Red(t_i)^\Sigma(r\leftarrow
u)\]
for $i=1,2$ and  $r,u\in W_\tau(X)$. Hence
\[\Sigma\models Red(t_1^\Sigma(r\leftarrow u))Red(t_2^\Sigma(r\leftarrow u))
\approx Red(t_1)^\Sigma(r\leftarrow u)Red(t_2)^\Sigma(r\leftarrow u).\]
  Thus we have
\[\Sigma\models Red(t_1^\Sigma(r\leftarrow u))Red(t_2^\Sigma(r\leftarrow
u))\approx
Red((t_1^\Sigma(r\leftarrow u)t_2^\Sigma(r\leftarrow u)))\]
\[\approx
Red(t^\Sigma(r\leftarrow u)).\]
Let us assume that $\Sigma\not\models t_1\approx t$ and $\Sigma\not\models t_2\approx
t$. Then  
\[\Sigma\models Red(t_1)^\Sigma(r\leftarrow u)Red(t_2)^\Sigma(r\leftarrow
u)\approx (Red(t_1)Red(t_2))^\Sigma(r\leftarrow u)\]
\[=Red(t)^\Sigma(r\leftarrow
u),\]
which proves (\ref{eq.4}) in this case.

Assume that $\Sigma \models t_1\approx t$. Then $\{1,2\}\subseteq
PEss(t,\Sigma)$ implies
 $\Sigma \models t_1\approx t_2$. Then
we have $\Sigma\models Red(t_1)\approx Red(t_2)$ and  
$\Sigma\models Red(t_1t_2)\approx Red(t_1)Red(t_2)\approx Red(t_1).$ Hence
\[\Sigma\models Red(t_1)^\Sigma(r\leftarrow u)Red(t_2)^\Sigma(r\leftarrow u) 
\] \[ \approx Red(t_1)^\Sigma(r\leftarrow u)\approx Red(t)^\Sigma(r\leftarrow u).\] 
\end{proof}

\begin{theorem}\label{t5.1}
The variety $IG=Mod(\{x_1x_1\approx x_1\})$ of idempotent groupoids is
stable.
\end{theorem}

\begin{proof}
We put $\Sigma =Id(Mod(\{x_1x_1\approx x_1\})$.  We have to prove  
(\ref{eq.2})
when $\Sigma\models t\approx s$ and  $r\in SEss(t,\Sigma)\cap SEss(s,\Sigma)$. Without loss of generality, let us assume that $Depth(t)\leq Depth(s)$.
 We shall proceed by induction on $Depth(t)$.

 Our inductive basis is $Depth(t)\leq 1$. Then clearly $t=s$ and (\ref{eq.2}) is
satisfied.

Assume that (\ref{eq.2}) is satisfied when $Depth(t)<k$ for some natural number
$k\geq 2$.

Let $Depth(t)=k$. Then $t=t_1t_2$ and $s=s_1s_2$. Lemma
\ref{l5.2} allows us to think that terms $t$ and $s$ are presented in their normal forms under $\rightarrow_{R}$, i.e. $t=Red(t)$ and $s=Red(s)$. Hence $\mbox{\small 1,2}\in PEss(t,\Sigma)\cap
PEss(s,\Sigma)$ and $\Sigma \not\models t_i\approx t$  and $\Sigma \not\models
s_i\approx s$ for $i=1,2$. This show that $\Sigma \not\models t_1\approx t_2$ 
and $\Sigma \not\models s_1\approx s_2$. Hence $\Sigma \models t_i\approx s_i$
for $i=1,2$. Now $Depth(t_i)<k$, our inductive assumption and Lemma  \ref{l5.2}
prove 
(\ref{eq.2}).
\end{proof}

\begin{theorem}\label{t5.2}
 The variety $CG=Mod(\{x_1x_2\approx x_2x_1\})$ of all the commutative
groupoids is stable.
\end{theorem}
\begin{proof}
 Let $\Sigma=\{x_1x_2\approx x_2x_1\}$. Let us note that
$\Sigma\models u\approx v$ implies $Len(u)=Len(v)$, $Depth(u)=Depth(v)$ and
$|Pos(u)|=|Pos(v)|$ for all $v,u\in W_\tau(X)$.

We shall prove (\ref{eq.2})  by induction on the depth of terms $t$ and $s$. 

Let  $Depth(t)=Depth(s)=0$. Then $t=s=x_1$ for some variable $x_1\in X$ and  (\ref{eq.2}) is obvious. 

Assume that (\ref{eq.2}) is satisfied when  $Depth(t)=Depth(s)<k$ for some
natural number $k$, $k>1$.

Let $Depth(t)=Depth(s)=k$,  $\Sigma\models t\approx s$ and  $r\in
SEss(t,\Sigma)\cap SEss(s,\Sigma)$. Let $n$ be a natural number such that
$t,s,r,u\in W_\tau(X_n)$ and let us denote by  $z_{n+1},\ldots,z_{n+p}\in W_\tau(X_n)$ 
all   subterms of $t,s$ or $r$ with depths equal to $1$ which are distinguished by
$\Sigma$, i.e. $\Sigma\not\models z_{n+i}\approx z_{n+j}$ when $i\neq j$. Using
inductive composition we obtain three new terms, namely:
\[t'=t(z_{n+1}\leftarrow x_{n+1},\ldots,z_{n+p}\leftarrow x_{n+p}),\]
\[s'=s(z_{n+1}\leftarrow x_{n+1},\ldots,z_{n+p}\leftarrow x_{n+p})\] and
\[r'=r(z_{n+1}\leftarrow x_{n+1},\ldots,z_{n+p}\leftarrow x_{n+p}).\]
Thus we have $t',s',r'\in W_\tau(X_{n+p})$ and $Depth(t')=Depth(s')=k-1<k$. It
is easy to see that $\Sigma\models t'\approx s'$ and  $r'\in SEss(t',\Sigma)\cap
SEss(s',\Sigma)$. Our inductive assumption implies 
$\Sigma\models t'^\Sigma(r'\leftarrow u)\approx s'^\Sigma(r'\leftarrow
u).$ Let us put $t''=t'^\Sigma(r'\leftarrow u)$ and
$s''=s'^\Sigma(r'\leftarrow
u)$. Then from $D_5$ it follows that
\[\Sigma\models t''(x_{n+1}\leftarrow z_{n+1},\ldots, x_{n+p}\leftarrow
z_{n+p})\] \[\approx 
s''(x_{n+1}\leftarrow z_{n+1},\ldots, x_{n+p}\leftarrow z_{n+p}).\]
Now, the equations
\[t^\Sigma(r\leftarrow u)=t''(x_{n+1}\leftarrow z_{n+1},\ldots, x_{n+p}\leftarrow
z_{n+p})\]
and 
\[s^\Sigma(r\leftarrow u)=s''(x_{n+1}\leftarrow z_{n+1},\ldots, x_{n+p}\leftarrow
z_{n+p})\]
complete the proof.
\end{proof}

\begin{remark}~~
$(i)$ It is  surprising  that the variety $CG$ of all commutative groupoids is stable, but the analogous variety of commutative semigroups is not stable, as shown by Theorem \ref{t3.1}. Hence the stability is not inherited by subvarieties of groupoids.

$(ii)$ Theorem \ref{t3.2} and the description of the lattice of the varieties of semigroups given in \cite{evans} show that if a variety $\mathcal V$  of semigroups is stable then all subvarieties of $\mathcal V$ are stable.
\end{remark}

Next we consider the following  varieties of groupoids:
\[V^{ijk}_{lm}=Mod(\{(x_ix_j)x_k\approx x_lx_m\})\quad and\quad W^{ijk}_{lm}=Mod(\{x_i(x_jx_k)\approx x_lx_m\}),\]
where $i,j,k,l,m\in\{1,2,3\}$.
\begin{theorem}\label{t5.3}
 The varieties of groupoids 
 $V^{ijk}_{lm}$ and 
$W^{ijk}_{lm}$ for $i,j,k,l,m\in\{1,2\}$ are
stable.
\end{theorem}
\begin{proof}
Since $Id(V^{ijk}_{lm})$ and $Id(W^{ijk}_{lm})$ are fully invariant congruences, they satisfy rules $D_1 - D_5$. Thus we have to prove that $\Sigma R_1$ is  satisfied in $Id(V^{ijk}_{lm})$ and $Id(W^{ijk}_{lm})$, i.e. that 
(\ref{eq.2}) is satisfied in $V^{ijk}_{lm}$ and $W^{ijk}_{lm}$.

Let $t,s,r\in W_\tau(X)$ be three terms for which $\Sigma\models t\approx s$ and $r\in SEss(t,\Sigma)\cap SEss(s,\Sigma)$. Suppose with no loss of generality  that $Depth(t)\leq Depth(s)$. 

If $\Sigma\models t\approx r$ then (\ref{eq.2}) is obvious. Thus we assume that  $\Sigma\not\models t\approx r$.  
\vspace{.5cm}

\noindent
{\bf Claim 1.} The varieties $V^{ijk}_{11}$ and $W^{ijk}_{11}$ for $i,j,k\in \{1,2\}$ are stable. 
\vspace{.5cm}

\noindent
In \cite{sht1} it is proved that $V^{121}_{11}$ is stable (see Proposition 3.1 of \cite{sht1}). In a similar way, one can prove that $V^{ijk}_{11}$ and $W^{ijk}_{11}$ for $i,j,k\in \{1,2\}$ are stable. 
\vspace{.5cm}

\noindent
{\bf Claim 2.}  The varieties $V^{ijk}_{12}$ and $W^{ijk}_{12}$ for $i,j,k\in \{1,2\}$ are stable.
\vspace{.5cm}

\noindent
We shall show that $V^{121}_{12}$ is stable by
induction on $Depth(t)$.

If $Depth(t)=0$ then  (\ref{eq.2}) is clearly satisfied.
Let $Depth(t)=1$. Then, with no loss of generality, we can assume that $t=x_1x_1$ or $t=x_1x_2$. 
Hence $\Sigma\models r\approx x_1$ or $\Sigma\models r\approx x_2$. Then  (\ref{eq.2}) follows from $D_5$.

Let us assume that (\ref{eq.2}) is satisfied when $Depth(t)<k$ for some natural
number $k$, $k>2$.

 Let $Depth(t)=k$. Then we  
have $t=t_1t_2$ with $1\leq Depth(t_i)< k$ for $i=1,2$. 
Then $\Sigma\models t\approx s$ implies that $s=s_1s_2$ or $s=(s_1s_2)s_1$ with
$\Sigma\models t_1\approx s_1$ and $\Sigma\models t_2\approx s_2$. Since $\Sigma\not\models t\approx r$ it follows that $r\in SEss(t_1,\Sigma)\cap SEss(s_1,\Sigma)$ or $r\in SEss(t_2,\Sigma)\cap SEss(s_2,\Sigma)$. Thus we have 
\[t^\Sigma(r\leftarrow u)=t_1^\Sigma(r\leftarrow u)t_2^\Sigma(r\leftarrow u),\]
\[s^\Sigma(r\leftarrow u)=s_1^\Sigma(r\leftarrow u)s_2^\Sigma(r\leftarrow u)\]
or
\[s^\Sigma(r\leftarrow u)=(s_1^\Sigma(r\leftarrow u)s_2^\Sigma(r\leftarrow u))s_1^\Sigma(r\leftarrow u),\] which proves (\ref{eq.2}), according to our inductive assumption. 

In a similar way one can show that $V^{211}_{12}$ and $V^{112}_{12}$ are stable varieties. By dual arguments we obtain that $W^{121}_{12}$, $W^{211}_{12}$ and $W^{112}_{12}$ are stable varieties.
\vspace{.5cm}

\noindent
{\bf Claim 3.} The varieties  $V^{111}_{ii}$ and $W^{111}_{ii}$ for $i\in\{1,2\}$ are stable.
\vspace{.5cm}

\noindent
We shall prove that $V^{111}_{11}$ is stable  by
induction on $Depth(t)$. 

If $Depth(t)=0$ or $Depth(t)=1$ then (\ref{eq.2}) can be proved as in the previous case.

Let us assume that (\ref{eq.2}) is satisfied when $Depth(t)<k$ for some natural
number $k$, $k>2$.

 Let $Depth(t)=k$. Then we  
have $t=t_1t_1$ or $t=t_1t_2$ with $1\leq Depth(t_i)< k$ for $i=1,2$. 

If $t=t_1t_1$ then $\Sigma\models t\approx s$ implies that $s=s_1s_1$ or $s=(s_1s_1)s_1$ with
$\Sigma\models t_1\approx s_1$. Since $\Sigma\not\models t\approx r$ it follows that $r\in SEss(t_1,\Sigma)\cap SEss(s_1,\Sigma)$.  Thus we have 
\[t^\Sigma(r\leftarrow u)=t_1^\Sigma(r\leftarrow u)t_1^\Sigma(r\leftarrow u),\]
\[s^\Sigma(r\leftarrow u)=s_1^\Sigma(r\leftarrow u)s_1^\Sigma(r\leftarrow u)\]
or
\[s^\Sigma(r\leftarrow u)=(s_1^\Sigma(r\leftarrow u)s_1^\Sigma(r\leftarrow u))s_1^\Sigma(r\leftarrow u),\] which proves (\ref{eq.2}), according to our inductive assumption. 

If  $t=t_1t_2$ then $\Sigma\models t\approx s$ implies that $s=s_1s_2$  with
$\Sigma\models t_1\approx s_1$ and $\Sigma\models t_2\approx s_2$. Since $\Sigma\not\models t\approx r$ it follows that $r\in SEss(t_1,\Sigma)\cap SEss(s_1,\Sigma)$ or $r\in SEss(t_2,\Sigma)\cap SEss(s_2,\Sigma)$. Thus we have 
\[t^\Sigma(r\leftarrow u)=t_1^\Sigma(r\leftarrow u)t_2^\Sigma(r\leftarrow u),\]
\[s^\Sigma(r\leftarrow u)=s_1^\Sigma(r\leftarrow u)s_2^\Sigma(r\leftarrow u)\]
 which proves (\ref{eq.2}) again, according to our inductive assumption. 
 
The varieties $V^{111}_{22}$ and $W^{111}_{22}$ are clearly stable.
\end{proof}

\section{S-stable varieties}\label{sec6}

Let us go back to the identities   (\ref{eq.1}). These identities guarantee 
stability of a variety of semigroups, satisfying one of them.  It is natural to expect that the identities 
(\ref{eq.1}) will provide for stability of the variety of groupoids. The next
proposition is a counterexample to that expectation. 

\begin{proposition}\label{p6.1}
 The varieties 
 $V_{lm}^{ijk}$ and  $W_{lm}^{ijk}$
 are not stable when $\{i,j,k\}=\{1,2,3\}$,\   $ l,m\in\{1,2,3\}$ and $l\neq
m$.
\end{proposition}
\begin{proof} Without loss of generality we prove that  $V_{23}^{123}=Mod(\Sigma)$ is not stable, where  $\Sigma=\{(x_1x_2)x_3\approx x_2x_3\}$. Let us put
$t=(x_3(x_1x_2))(x_2(x_1x_2))$,
$s=x_2(x_2(x_1x_2))$,
 $r=x_1x_2$ and $u=x_4$. Clearly, $\Sigma\models t\approx
s$. Since $PFic(t,\Sigma)=\{\mbox{\small 11,121}\}$ and  $PFic(s,\Sigma)=\emptyset$,  it follows that $P_r^t=\{\mbox{\small 12,22}\}$ and
$P_r^s=\{22\}$. Thus we have 
$t^{\Sigma}(r\leftarrow u)=(x_3x_4)(x_2x_3)$ and $s^{\Sigma}(r\leftarrow
u)=x_2(x_2x_4)$. Clearly,
$\Sigma\not\models t^{\Sigma}(r\leftarrow u)\approx s^{\Sigma}(r\leftarrow u)$.
Hence the variety $Mod(\Sigma)$ is not stable.
\end{proof}

Our aim in this
section is to define additional sufficient  conditions for stability  
such that if a variety of groupoids satisfies an identity among   (\ref{eq.1})
then it is  stable under these conditions. Also, we expect  the varieties $V_{11}$, $V_{22}$ and $V_{33}$ of semigroups to be included in this new concept of stability.
We are going to  define s-stable variety for an arbitrary type $\tau$.

 Let $t$ and $r$ be two terms and let 
 $EP_r^t=\{p\in P_r^t \mid   p\preceq q\in Pos(t)  \Rightarrow  q\in
PEss(t,\Sigma)  \}$   be the set of  all the
minimal elements in $\Sigma P_r^t$ whose successors are $\Sigma$-essential in $t$.

\begin{definition}%
\label{d6.1} Let $r,s,t\in W_\tau(X)$ be  terms of
type $\tau$. The
\emph{essential composition} of the terms $t$ and $r$ by $s$ is defined as follows
 
\begin{enumerate}
\item[(i)] $t(r* s)=t$ if $EP_r^t=\emptyset$;

\item[(ii)] $t(r* s)=s$ if $\Sigma\models t\approx r$,  and

\item[(iii)] $t(r* s)=f(t_1(r*
s),\ldots,t_n(r* s))$, if
$t=f(t_1,\ldots,t_n)$ and $\Sigma\not\models t\approx r$.
\end{enumerate}
\end{definition}

\begin{example}\label{ex2}
 Let us consider the terms $t=(x_3(x_1x_2))(x_2(x_1x_2))$ and  $r=x_1x_2$ from the Proposition \ref{p6.1} and let
$s=x_4$. 
Then we have $EP_r^t=\{\mbox{\small 22}\}$ and $t(r*s)=(x_3(x_1x_2))(x_2x_4).$ On the other side $P_r^t=\{\mbox{\small 12,22}\}$ implies $t^\Sigma(r\leftarrow s)=(x_3x_4)(x_2x_4)$. Clearly, 
\[\Sigma\not\models t(r*s)\approx t^\Sigma(r\leftarrow s).\]
\end{example}

\begin{definition}\label{d6.2}
A set   $\Sigma$ of identities  is $SR$-deductively closed if
it satisfies the rules $D_1,D_2,D_3,D_5$ and
\begin{enumerate}

\item[$SR_1$] \emph{(Star\   Replacement)}
\[{\left(\begin{array}{c}
r,t,s,u\in W_\tau(X)\ \&\ (t\approx s\in\Sigma)\ \&\   \\
                 (EP_r^t\neq\emptyset)\ \&\  (EP_r^s\neq\emptyset)
\end{array}\right)}\Rightarrow  t(r* u)\approx
s(r*u)\in\Sigma.\]
\end{enumerate}
\end{definition}

For any set of identities $\Sigma$ the smallest $SR$-deductively closed set
containing $\Sigma$ is called the \emph{ $S
R$-closure} of $\Sigma$, and denoted by $S R(\Sigma).$
 For $t\approx s\in Id(\tau)$
we say $\Sigma\vdash_{S R} t\approx s $ ($``\Sigma$ $S R
$-proves $t\approx s"$) if there is a sequence of identities
$t_1\approx s_1,\ldots,t_n\approx s_n$, such that each identity
belongs to $\Sigma$ or is a result of applying any of the derivation
rules
   $D_1, D_2, D_3,  D_5$ or $S R_1$
  to previous identities
in the sequence and the last identity $t_n\approx s_n$ is
$t\approx s.$

Let $t\approx s$ be an identity and $\mathcal{A}$ be an algebra of
type $\tau$. $\mathcal{A}\models_{S R} t\approx s$ means that
$\mathcal{A}\models t(r* v) \approx
s(r* v)$ for every $r\in SEss(t,\Sigma)\cap
SEss(s,\Sigma)$ and $v\in W_\tau(X)$.
 For $t,s\in W_\tau(X)$ we say
$\Sigma\models_{S R} t\approx s$ (read: ``$\Sigma$ $S
R$-yields $t\approx s$") if, given any algebra $\mathcal{A}$,
$\mathcal{A}\models_{S R} \Sigma\quad\Rightarrow\quad \mathcal{A}\models_{S R}
t\approx s.$

As in \cite{sht1} (see Theorem 3.4 and Theorem 3.6) one can prove that $SR$  is a
closure operator, and prove a completeness theorem that 
$\Sigma\models_{SR}t\approx s \iff \Sigma\vdash_{SR}t\approx s$. 

\begin{theorem}\label{t6.1} For each set of identities $\Sigma$ the closure
 $S R(\Sigma)$ is a fully invariant congruence.
 \end{theorem}
\begin{proof}\ \ It is enough to  prove that $SR(\Sigma)$  satisfies the rule
 $D_4$. Let $r\in W_\tau(X)$, $t\approx s
\in\Sigma$ and $p\in Pos(r)$.
If $p\notin PEss(r,\Sigma)$, then  we have $r(p;v)\approx r(p;w)\in SR(\Sigma)$
for
all terms $v,w\in W_\tau(X)$.
Let $p\in PEss(r,\Sigma)$ and let $n$ be a natural number such that
$r,t,s\in W_\tau(X_n)$. Denote by  $v=r(p;x_{n+1})$ and $u=x_{n+1}$.  Clearly,
$u\in Sub(v)$ and $EP_u^v=\{p\}$.
We have
$v(u*t)=r(p;t)$ and $v(u*s)=r(p;s)$.
Now from $S R_1$ we obtain
$v(u*t)\approx v(u*s)\in SR(\Sigma),$   i.e. $r(p;t)\approx
r(p;s)\in SR(\Sigma).$
 \end{proof}

Since $EP_r^t\subseteq P_r^t$ and $EP_v^s\subseteq P_v^s$ we have 
$t\approx s\in SR(\Sigma)\ \Rightarrow\ t\approx s\in \Sigma R(\Sigma),$ 
for each identity $t\approx s\in Id(\tau)$. Thus we obtain the following
inclusions
$D(\Sigma)\subseteq SR(\Sigma)\subseteq \Sigma R(\Sigma)$ for each
$\Sigma\subseteq Id(\tau)$. Hence each stable variety is s-stable one.

  \begin{definition}\label{d6.3}
A set of identities $\Sigma$ is called an \emph{ s-globally
invariant
congruence} if it is $SR$-deductively closed.

A variety $V$ of type $\tau$ is called \emph{s-stable} if
$Id(V)$   is an s-globally
invariant congruence.
\end{definition}

\begin{proposition}\label{p6.2}
There exist sets $\Sigma_1$ and $\Sigma_2$ of identities  such that $D(\Sigma_1)\subsetneqq SR(\Sigma_1)$ and 
$SR(\Sigma_2)\subsetneqq \Sigma R(\Sigma_2)$.
\end{proposition}
\begin{proof}
First, let   $\Sigma_1=
\{x_1(x_2x_3) \approx (x_1x_2)x_3 \}$ be the set of identities which
define the 
variety $SG=Mod(\Sigma_1)$ of semigroups.
Clearly,
$Id(SG)=D(\Sigma_1).$
 Let us set  $t=((x_1x_2)x_1)x_2$, $ s=(x_1x_2)(x_1x_2)$,
$r=x_1x_2$ and $u=x_3$.  Clearly $\Sigma_1\models t\approx s$.   Since
$EP_r^t=\{\mbox{\small 11}\}$ and $EP_r^s=\{\mbox{\small 1,2}\}$,
  we obtain
$t(r* u)=(x_3x_1)x_2$ and  $s(r*u)=x_3x_3.$ 
 Hence
$\Sigma_1\not\models  t(r* u)\approx s(r*u)$. 
Consequently,
$D(\Sigma_1)$ is a proper subset of $SR(\Sigma_1)$ and $Mod(SR(\Sigma_1))$ is a proper
subvariety of
$SG$.

Second, let $\Sigma_2=\{(x_1x_2)x_3\approx x_2x_3\}$. Let us consider
the terms $t$,
$s$ and
 $r$  considered in Proposition \ref{p6.1}. It is easy to see that 
$EP_r^t=EP_r^s=\{\mbox{\small 22}\}$. Thus we have 
$ t(r* u)=(x_3(x_1x_2))(x_2x_3)$, $s(r*u)=x_2(x_2x_3)$ and hence
 $\Sigma_2\models t(r* u)\approx s(r*u)$, but 
$\Sigma_2\not\models t^{\Sigma_2}(r\leftarrow u)\approx s^{\Sigma_2}(r\leftarrow u)$.
\end{proof}

\begin{lemma}\label{l6.1}
 Let $x_i\in X$ be a $\Sigma$-essential variable which occurs once in the
term $t\in W_\tau(X)$. Then the variable $x_i$ is $\Sigma$-essential in
$Red(t)$ with a unique occurrence.
\end{lemma}
\begin{proof} According to Theorem \ref{t4.1}, it is enough to prove that
$x_i\in X$ is $\Sigma$-essential in $r$ with unique occurrence when
$t\rightarrow_{R} r$.  Corollary 3.8 of \cite{sht} and 
Corollary \ref{l4.1}    imply $x_i\in Ess(r,\Sigma)$.

Let $t\rightarrow_R r$, $r=t(p;u)$, $s=sub_t(p)$ and $\Sigma\models\ s\approx u$ where $u$ is $\Sigma$-minimal. Let $q$ be the unique position on which $x_i$ occurs in $t$. Since $q$ is a position of a variable it follows that $q\not\prec p$.

If $p\prec q$ then the unique occurrence of $x_i$ in $r$ follows by the $\Sigma$-minimality of $u$. If $p\not\prec q$ then $x_i$ occurs once on the position $q\in Pos(r)$ in $r$.
\end{proof}

\begin{lemma}\label{l6.2} If $\Sigma=\{ x_1(x_2x_3)\approx x_ix_j\}$  with $1\leq i\leq j\leq 3$ then
 \[\Sigma\models t(r*u)\approx Red(t)^\Sigma(r\leftarrow u)\] for all $t,r,u\in
W_\tau(X)$. 
\end{lemma}
\begin{proof}
 Let $n$ be a natural number such that $r,t,u\in W_\tau(X_n)$. If
$EP_r^t=\emptyset$ then $P_r^t=\emptyset$ and we are done.
 Let
  $EP_r^t=\{p_1,\ldots,p_m\}$ and let us put $s=t(p_1,\ldots,p_m;
x_{n+1}\ldots x_{n+m})$. Clearly $ x_{n+1}\ldots x_{n+m}\in Ess(s,\Sigma)$ and
$x_{n+i}$ occurs only once in $s$ for $i=1,\ldots,m$. From Lemma \ref{l6.1} it follows 
that $ x_{n+1}\ldots x_{n+m}\in Ess(Red(s),\Sigma)$ and $x_{n+i}$ occurs only
once in $Red(s)$ for $i=1,\ldots,m$. If we suppose that there is a term $v$ such
that $\Sigma \models r\approx v$ and $v\in Sub(Red(s))$ then there is $w\in
Sub(s)$ such that $\Sigma\models v\approx w$. Since $Ess(v,\Sigma)\subseteq 
Ess(Red(s),\Sigma)$ it follows that $Ess(v,\Sigma)\subseteq  Ess(s,\Sigma)$. Then
from Theorem 2.13 of \cite{sht1} it follows that $v\in EP_r^s\subseteq EP_r^t$
which is a contradiction. Hence $\Sigma\not\models r\approx v$ for all $v\in
Sub(Red(s))$. Consequently $EP_r^s=P_r^s=\emptyset$ and we obtain
$ t(r*u)=s(x_{n+1}\leftarrow u,\ldots, x_{n+m}\leftarrow u)$ and
\[Red(s)^\Sigma(r\leftarrow u)=Red(s)(x_{n+1}\leftarrow u,\ldots, x_{n+m}\leftarrow
u).\]
From Corollary \ref{l4.1} we have 
\[\Sigma\models 
s(x_{n+1}\leftarrow
u,\ldots, x_{n+m}\leftarrow u)\approx Red(s)(x_{n+1}\leftarrow
u,\ldots, x_{n+m}\leftarrow u)\]
which completes the proof.
\end{proof}

\begin{lemma}\label{l6.3}
If $\Sigma=\{(x_1x_2)x_3\approx x_ix_j\}$  with $1\leq i\leq j\leq 3$  then  the normal form under the reduction
$\rightarrow_{R}$ of a term $t\in W_\tau(X)$ is presented in the following form:
\begin{equation}\label{eq.5} Red(t)=x_{i_1}(x_{i_2}(\ldots (x_{i_{n-1}}x_{i_n})\ldots
)),\end{equation} where   $x_{i_m}\in var(t)$ for $m=1,\ldots,n$.
\end{lemma}
\begin{proof} Let   $\mathcal V$ be the variety defined by $\Sigma$, i.e.
$\mathcal V=Mod(\Sigma)$.
 We shall prove the lemma when $\Sigma\models (x_1x_2)x_3\approx
x_1x_2$. The other cases follow by dual arguments. 

So, let us   consider the term  $t=(x_1x_2)x_3$. Then we have $2\notin
PEss(t,\Sigma)$. Hence $Red(t)=x_1x_2$ and we are done.

 Assume that if $Depth(t)<k$, for some natural number $k, k>2$ then $Red(t)$ is
presented in the form of (\ref{eq.5}).

Let $Depth(t)=k$. Then we have $t=t_1t_2$ with $t_1,t_2\in W_\tau(X)$ and
$0\leq Depth(t_i)<k$ for $i=1,2$. Clearly  $\Sigma\models Red(t)\approx
Red(t_1)Red(t_2)$. From the inductive assumption we know that $Red(t_1)$ and
$Red(t_2)$ are presented in the form of (\ref{eq.5}).  If $Red(t_1)=x_{i_1}$ then we are done. Let $Depth(Red(t_1))\geq 1$ and $Red(t_1)=x_{i_1}t_{12}$ for some $t_{12}\in W_\tau(X)$. Then   
\[\Sigma \models Red(t)\approx  (x_{i_1}t_{12})Red(t_2)\approx 
x_{i_1}t_{12}=Red(t_1)\] 
which completes the proof.
\end{proof}

By dual arguments one can prove the following lemma.

\begin{lemma}\label{l6.4}
If $\Sigma=\{ x_1(x_2x_3)\approx x_ix_j\}$  with $1\leq i\leq j\leq 3$  then  the normal form under the reduction
$\rightarrow_{R}$ of a term $t\in W_\tau(X)$ is presented in the following form:
\[ Red(t)=(\ldots ((x_{i_1}x_{i_2})x_{i_3})\ldots)x_{i_n},\]
 where   $x_{i_m}\in var(t)$ for $m=1,\ldots,n$.
\end{lemma}

\begin{theorem}\label{t8}
The varieties of semigroups $V_{11}$, $V_{22}$ and $V_{33}$ are s-stable (see Proposition \ref{p1}).
\end{theorem}
\begin{proof} We shall prove that $V_{11}$ is an s-stable variety. To show that $\Sigma=Id(V_{11})$ is $SR$-deductively closed, i.e. $SR(\Sigma)=\Sigma$, we let $r,s,t$ be three terms such that $t\approx s\in\Sigma$, $EP_r^t\neq\emptyset$ and $EP_r^s\neq\emptyset$. 
 We have to prove 
 \begin{equation}\label{eq5}
 \Sigma\models t(r*u)\approx s(r*u).
 \end{equation}
 
 If $Depth(t)\leq 1$ then we have 
 \[\Sigma\models t\approx s\Longrightarrow t=s\] and  (\ref{eq5}) is obviously satisfied.
 
 Let $Depth(t)\geq 2$ and $Depth(s)\geq 2$.
 Since $x_1x_2x_3\approx x_1x_1\in\Sigma$, the set of $\Sigma$-essential positions in each term $w$ consists of all strings over $\{1\}$ which belong to $Pos(w)$, including the empty string $\varepsilon$. Consequently, for each term $r$ we have $EP_r^w=\emptyset$ or $EP_r^w=\{p_w\}$, where $p_w$ is the longest string over $\{1\}$ in $Pos(w)$. 
 
 Next,  $EP_r^t\neq\emptyset$ and $EP_r^s\neq\emptyset$ imply  $EP_r^t=\{p_t\}$ and $EP_r^s=\{p_s\}$. Since $p_t$ and $p_s$ are the longest strings in $Pos(t)$ and $Pos(s)$, respectively it follows that $r$ is a variable and $r=first(t)=first(s)$. Thus,  (\ref{eq5}) follows by $D_5$.

In a similar way one can prove that $V_{33}$ is an s-stable variety. The proof that  $V_{22}$ is a stable variety is left to the reader.
\end{proof}

\begin{theorem}\label{t6.2}
The varieties of groupoids
 $V^{ijk}_{lm}$ and 
$W^{ijk}_{lm}$ for $i,j,k,l,m\in\{1,2,3\}$ are
s-stable.

 \end{theorem}
\begin{proof}
If $i,j,k,l,m\in\{1,2\}$ we are done because of Theorem \ref{t5.3}.
\vspace{.5cm}

\noindent
{\bf Claim 1.} $V^{123}_{lm}$ and 
$W^{123}_{lm}$  with $1\leq l\leq m\leq 3$  are s-stable varieties.
\vspace{.5cm}

\noindent
We are going to   prove that $V_{12}^{123}=Mod(\Sigma)$ is s-stable, where  $\Sigma=\{(x_1x_2)x_3\approx
x_1x_2\}$. 
Lemma \ref{l6.2} implies $t(r*u)=Red(t)^\Sigma(Red(r)\leftarrow u)$
and it is enough to prove 
\begin{equation}\label{eq.6}
 \Sigma\models Red(t)^\Sigma(Red(r)\leftarrow u)\approx
Red(s)^\Sigma(Red(r)\leftarrow u)
\end{equation}
when $\Sigma\models t\approx s$, $r\in SEss(t,\Sigma)\cap SEss(s,\Sigma)$ and $u\in W_\tau(X)$.

Let $t,s,r\in W_\tau(X)$ be three terms for which $\Sigma\models t\approx s$ and $r\in SEss(t,\Sigma)\cap SEss(s,\Sigma)$. Suppose with no loss of generality  that $Depth(t)\leq Depth(s)$.

We argue by induction on  $Depth(t)$. If $\Sigma\models t\approx r$ then (\ref{eq.6}) is obvious.

Assume that $\Sigma\not \models t\approx r$.

Let $Depth(t)=1$. Then, without loss of generality we can assume that $t=x_1x_2$. Hence $\Sigma\models r\approx x_1$ or $\Sigma\models r\approx x_2$, and (\ref{eq.6}) follows from $D_5$.

Assume that  for some natural number  $k\geq 2$, if  $Depth(t)<k$ then  (\ref{eq.6}) is
satisfied.

Let $Depth(t)=k$.  From Lemma \ref{l6.3} it follows that
\[Red(t)=x_{i_1}(x_{i_2}(\ldots (x_{i_{n-1}}x_{i_n})\ldots ))\]
and 
 \[Red(s)=x_{j_1}(x_{j_2}(\ldots (x_{j_{m-1}}x_{j_m})\ldots )),\]
 where $x_{i_l}\in var(t)$ and $x_{j_k}\in var(s)$ for $l=1,\ldots,n$ and $k=1,\ldots,m$.
 Clearly $x_{i_1}=x_{j_1}$ because
$\Sigma\models t\approx s$ and $\mbox{\small 1}\in PEss(t,\Sigma)\cap PEss(s,\Sigma)$.

 If $Red(r)=x_{i_1}$ then  we are done because of $D_5$.

If  $Red(r)\neq x_{i_1}$ then $r\in SEss(t_2,\Sigma)\cap SEss(s_2,\Sigma)$ where 
\[t_2=x_{i_2}(\ldots (x_{i_{n-1}}x_{i_n})\ldots )\quad and\quad s_2=x_{j_2}(\ldots (x_{j_{m-1}}x_{j_m})\ldots ).\]

Clearly $\Sigma\models t_2\approx s_2$ and we have 
\[Red(t)^\Sigma(Red(r)\leftarrow u)= x_{i_1}Red(t_2)^\Sigma(Red(r)\leftarrow u),\]
and 
\[Red(s)^\Sigma(Red(r)\leftarrow u)= x_{i_1}Red(s_2)^\Sigma(Red(r)\leftarrow u)\]
for each $u\in W_\tau(X)$,
which  together with our inductive assumption prove  (\ref{eq.6}).

\vspace{.5cm}

\noindent
{\bf Claim 2.} $V^{123}_{lm}$ and 
$W^{123}_{lm}$  with $1\leq m<l\leq 3$  are s-stable varieties.
\vspace{.5cm}

We shall show that $V^{123}_{31}$ is s-stable. Thus we have 
\[\Sigma\models  x_1x_3\approx (x_3(x_4x_5))x_1\approx (x_1x_2)(x_3(x_4x_5))\approx\] \[ ((x_2x_6)x_1)(x_3(x_4x_5))\approx (x_3(x_4x_5))(x_2x_6)\approx (x_2x_6)x_3\approx x_3x_2.\]
Hence $\Sigma\models x_1x_3\approx x_3x_2\approx x_2x_4$.  Hence, if $Depth(t)\geq 1$ then without loss of generality we can assume that $\Sigma\models t\approx Red(t)=x_1x_2$ with $PEss(Red(t))={\varepsilon}$. Consequently, for each term $r$ we have $EP_r^t=\emptyset.$
\end{proof}

\end{document}